\documentclass[12pt]{article}
\usepackage{amsmath}
\usepackage{amscd}
\usepackage{graphicx}
\usepackage{amssymb,amsfonts,amsthm,amscd,mathrsfs}
\usepackage{verbatim}
\usepackage[T2A]{fontenc}
\usepackage[utf8]{inputenc}
\usepackage[russian]{babel}
\usepackage{bm}
\usepackage{authblk}

\renewcommand{\b}{\beta}

\newcommand{\e}{\varepsilon}

\newcommand{\E}{\mathsf{E}}

\begin{document}

\title{On stochastic expansions of empirical distribution function of residuals in autoregression schemes
}

\author{M.~V.~Boldin    
\footnote{Moscow State Lomonosov Univ., Dept. of Mech. and Math., Moscow, Russia\\
e-mail: boldin$_{-}$m@hotmail.com}}

\date{ }
\maketitle

\textbf{Abstract}

We consider a stationary linear AR($p$) model with unknown mean. The autoregression parameters as well as the distribution function (d.f.) $G$ of innovations are unknown. The observations contain gross errors  (outliers).  The distribution of outliers is unknown and arbitrary, their intensity is
 $\gamma n^{-1/2}$ with an unknown $\gamma$, $n$ is the sample size. 
 
 The assential problem  in such situation is to test the normality of innovations.
 Normality, as is known, ensures the optimality properties of widely used least squares procedures.  To construct and study a Pearson chi-square type test for normality we estimate the unknown mean and the autoregression parameters. Then, using the estimates, we find the residuals in the autoregression. Based on them,  we construct a kind of empirical distribution function (r.e.d.f.) , which is a counterpart of the (inaccessible) e.d.f. of the autoregression innovations.
Our Pearson's satatistic is the functional from r.e.d.f. Its asymptotic distributions under the hypothesis and the local alternatives are determined by the asymptotic behavior of  r.e.d.f. 


In the present work, we find and substantiate in details the stochastic expansions of the r.e.d.f. in two situations. In the first one d.f. $ G (x) $ of innovations does not depend on $ n $. We need this result to investigate test statistic under the hypothesis. In the second situation $ G (x) $ depends on $ n $ and has the form of a mixture $ G (x) = A_n (x) = (1-n ^ {- 1/2}) G_0 (x) + n ^ { -1/2} H (x). $ We need this result to study the power of test under the local alternatives. 



\section{
Постановка задачи и результаты }
В этой работе мы рассматриваем стационарную AR($p$) модель с ненулевым средним
\begin{equation}
v_t = \b_1 v_{t-1} + \dots + \b_p v_{t-p} +\nu + {\e}_t, \quad  t \in \mathbb{Z} .
\end{equation}

В (1.1) $\{{\e}_t\}$ -- независимые одинаково распределенные случайные величины (н.о.р.сл.в.) с неизвестной функцией распределения (ф.р.) $G(x)$;
 $\E \e_1 = 0$,  $0<\E {\e}_1^2 < \infty$; $\bm{\b} = (\b_1, \dots, \b_p)^T \in  \mathbb{R}^p $-- вектор неизвестных параметров, таких что корни соответствующего (1.1) характеристического уравнения по модулю меньше единицы;  $\nu$ -- неизвестное среднее,  $\nu \in \mathbb R^1$.\\
 Эти условия дальше всегда предполагаются выполненными и особо не оговариваются.\\
 Мы рассматриваем модель (1.1) с выбросами в наблюдениях. А именно, предполагается, что наблюдаются величины
$$
y_t=v_t+z^{\gamma_n}_t {\xi}_t,\quad t  = 1-p, \dots, n, \eqno(1.2)
$$
где $v_{1-p}, \dots, v_n$ -- выборка из стационарного решения $\{v_t\}$ уравнения  (1.1); $\{z^{\gamma_n}_t\}$ н.о.р.сл.в., имеющие распределение Бернулли, т.е. принимающие значения 1 и 0, причем вероятность единицы
 $\gamma_n$,
 $$
 \gamma_n = \min(1, \frac{\gamma}{\sqrt{n}}) ,\quad  \gamma \ge 0 \text{\; неизвестно.}
 $$
Кроме того, $\{ \xi_t \}$ -- н.о.р.сл.в. с произвольным и неизвестным распределением $\Pi$. Последовательности $\{v_t\},\,\{z_t^{\gamma}\},\,\{\xi_t\}$ независимы между собой.

Переменные $\{ \xi_t \}$ интерпретируются как выбросы (засорения), $\gamma_n$ уровень засорения. Для $\gamma = 0$ мы получаем модель (1.1) без засорений.

Модель (1.2) -- локальный вариант хорошо известной модели засорения данных во временных рядах, см. \cite{MartYoh86}.\\
Перепишем уравнение (1.1) в удобном для дальнейшего рассмотрения виде. Для этого определим константу $\mu$ соотношеним
$$
\nu=(1-\b_1-\ldots-\b_p)\mu,
$$
тогда
$$
v_t -\mu= \b_1 (v_{t-1}-\mu) + \dots + \b_p( v_{t-p} -\mu)+ {\e}_t, \quad  t \in \mathbb{Z.}
$$

Если положить $u_t:=v_t-\mu,$ то
$$
v_t=\mu+u_t,\quad u_t=\b_1 u_{t-1} + \dots + \b_pu_{t-p}+ {\e}_t, \quad  t \in \mathbb{Z.}\eqno(1.3	)
$$
Последовательность $\{u_t\}$ в (1.3) -- авторегрессионная последовательность с нулевым средним и конечной дсперсией.\\
Построим по наблюдениям $\{y_t\}$ из (1.2) оценки ненаблюдаемых $\{\e_t\}$.\\ 
Далее $\Gamma$ -- любое конечноое неотрицательное число. Пусть $\hat\mu_n$ будет любая оценка $\mu$, для которой 
 последовательнсть
$$
n^{1/2}(\hat\mu_n-\mu)=O_P(1),\quad n\to\infty, \text{ равномерно по}\,\,\,\gamma \leq\Gamma.\eqno(1.4)
$$
Например, в качестве оценки $\hat\mu_n$ можно взять обычную М-оценку, построенную по $\{y_t\}$ так же, как строится М-оценка параметра сдвига по н.о.р. данным.
 
 Положим
$$
\hat u_t=y_t-\hat\mu_n,\quad t=1,\ldots,n.
$$
Пусть $\hat{\bm\b}_n= (\hat \b_{1n}, \dots, \hat \b_{pn})^T$ будет любая оценка $\bm\b$,
для которой 
 последовательность
 $$
 n^{1/2}(\hat{\bm\b}_n-\bm\b)=O_P(1),\quad n\to \infty, \text{ равномерно по}\,\,\gamma\leq \Gamma.\eqno(1.5)
 $$
Примером подходящих оценок $\hat{\bm\b}_n$ могут служить GM-оценки, построенные по $\{\hat u_t\}$ так же, как они строятся по $\{u_t\}$. Асимптотические свойства таких оценок аналогичны свойствам GM-оценок из Раздела  2.4 \cite{Bold.Petr.}. \\
Положим
$$
\hat \e_t =  \hat {u}_t -  \hat \b_{1n}\hat { u}_{t-1} - \dots - \hat \b_{pn}\hat{ u}_{t-p},\quad t = 1, \dots, n.\eqno(1.6)
$$
Величины $\{\hat \e_t\}$ называются остатками. 

Ведем остаточную эмпирическую функцию распределения ( о.э.ф.р.)
$$
\hat G_n(x) = n^{-1} \sum_{t=1}^n I(\hat \e_t \le x),\quad  x\in \mathbb{R}^1.
$$
 Здесь и в дальнейшем $I(\cdot)$ обозначает индикатор события.\\
Функция $\hat G_n(x)$ -- аналог гипотетической э.ф.р.
$$
G_n(x) = n^{-1} \sum_{t=1}^n I(\e_t \le x)
$$
ненаблюдаемых величин $\e_{1}, \dots, \e_n$.\\ 

Мы будем рассматривать две ситуации. В первой ф.р. $G(x)$ величин $\{\e_t\}$ не зависит от числа наблюдений $n$. Наша первая задача -- доказать следующую теорему 1.1.
\newtheorem{Th}{Теорема}[section]
\begin{Th}
{\sl 
Предположим, что ф.р. $G(x)$ имеет дифференцируемую производную $g(x)=G'(x)$, и $\sup_x |g'(x)|<\infty$. Пусть выполнены соотношения (1.4) -- (1.5). Тогда при любых фиксированных $x,\,\Gamma \ge 0$ и любом 
 $\delta > 0$
$$
\sup_{\gamma \leq \Gamma}  \Prob(|n^{1/2}[\hat G_n(x) -G_n(x)]- g(x)\delta (\beta) n^{1/2}(\hat\mu_n-\mu) -  \gamma\Delta (x, \Pi) | > \delta) \to 0,
\quad n \to \infty.
$$
Здесь сдвиг, определяемый засорениями,
$$
\Delta (x, \Pi) = \sum_{j=0}^p [\E G(x + \b_j \xi_1) - G(x)],\, \b_0 =-1,\,
$$
а $\delta (\beta)=1- \b_{1}-\ldots- \b_{p}$.
}
\end{Th}

Во второй ситуации относительно $G(x)$ предполагается, что она зависит от $n$, $G(x)=A_n(x)$, и $A_n(x)$ задается смесью
$$
A_n(x)=(1-n^{-1/2})G_0(x)+n^{-1/2}H(x).\eqno(1.7)
$$
В (1.7) функции распределения $G_0(x)$ и $H(x)$ удовлетворяют следующему условию.

{\bf Условие (А).}
Функции распределения $G_0(x)$ и $H(x)$ имеют средние ноль, конечные дисперсии. Они дважды дифференцируемы с ограниченными вторыми призводными.\\

Наша вторая задача -- доказать в ситуации два следующую теорему.
\begin{Th}
{\sl 
Предположим, что ф.р. $G(x)$ зависит от $n$ и задается соотношением (1.7). Пусть выполнено Условие (А), и $g_0(x)=G_0'(x)$. Пусть выполнены соотношения (1.4) -- (1.5). Тогда при любом фиксированном $x$ и любом 
 $\delta > 0$
$$
\sup_{\gamma\leq \Gamma}  \Prob(|n^{1/2}[\hat G_n(x) -G_n(x)]- g_0(x)\delta (\beta) n^{1/2}(\hat\mu_n-\mu) -  \gamma\Delta_0 (x, \Pi) | > \delta) \to 0,
\quad n \to \infty.\eqno(1.8)
$$
Здесь сдвиг
$$
\Delta_0 (x, \Pi) = \sum_{j=0}^p [\E G_0(x + \b_j \xi_1) - G_0(x)],\, \b_0 =-1,\,
$$
и, как в теореме 1.1, $\delta (\beta)=1- \b_{1}-\ldots- \b_{p}$.
}
\end{Th}

{\bf Замечание 1.1.}

В разложении (1.8) формально отсутствует ф.р. $H(x)$. Но $H(x)$ присутствует в (1.8) неявно, от нее зависит предельное распределение $G_n(x)$. А именно, обозначим через $D[0,1]$ пространство Скорохода функций без разрывов 2-го рода, определенных на $[0,1]$. Пусть $v(t),\,t\in [0,1]$, будет броуновский мост. Тогда, см. \cite{Chib65}, $
\hat v_n(t):=n^{1/2}[G_n(G_0^{-1}(t))-t]$ слабо сходится в $D[0,1]$ к $v(t)+[H(G_0^{-1}(t))-t],\,n\to \infty$.

{\bf Замечание 1.2.}

Теорема 1.2 была анонсирована без доказательства в \cite{Bold.Arx.2}. Цель настоящей работы -- дать полное доказательство теорем 1.1 и 1.2. Разумеется,  теорема 1.1 следует из теоремы 1.2 при $G_0(x)=H(x)=G(x)$, но удобно  сформулировать именно две теоремы, т.к. доказательство теоремы 1.1 менее громоздко и на нем удобно объяснить доказательство теоремы 1.2.\\
{\bf Замечание 1.3.}

В силу сформулированных теорем предельное распределение о.э.ф.р. $\hat G_n(x)$ зависит от распределения оценки $\hat\mu_n$. Предположим, что ф.р. $G(x)$ симметрична относительно нуля.
Будем брать в такой ситуации оценкой $G(x)$ симметризованную оценку
$$
\hat S_n(x):=\frac{\hat G_n(x)+1-\hat G_n(-x)}{2}.
$$
Положим
$$
\Delta_S(x, \Pi):=\frac{\Delta_0(x, \Pi)-\Delta_0(-x, \Pi)}{2}.
$$
Пусть
$$
 S_n(x):=\frac{ G_n(x)+1- G_n(-x)}{2}.
$$

Теоремы 1.1 и 1.2  прямо влекут
\newtheorem{Corollary}{Следствие
}[section]
\begin{Corollary}
{\sl
При условиях теоремы 1.1 для симметричной $G(x)$ или в условиях теоремы 1.2 при симметричной $G_0(x)$
$$
\sup_{\gamma \le \Gamma} \Prob(|n^{1/2}[\hat S_n(x) -S_n(x)]-  \gamma\Delta_S(x, \Pi) | > \delta) \to 0,
\quad n \to \infty.
$$
}
\end{Corollary}
Применения сформулированных теорем и следствия к построению и исследованию симметризованных тестов типа хи-квадрат Пирсона для гипотез относительно $G(x)$  Пирсона были анонсированы в \cite{Bold.Arx.2}, \cite{Bold.Arx.3}. В частности, там описаны тесты для проверки нормальности $G(x)$.

При этом, теорема 1.1 нужна, чтобы находить предельное распределение тестовых статистик при гипотезах, а теорема 1.2 -- при локальных альтернативах.

\section{Доказательство теоремы 1.1.}
Докажем теорему 1.1 для $p=1$, доказательство в общем случае отличается лишь громоздкими, но непринципиальными  техническими детатями. \\
Переобозначим $\b_1$ как $\b$, тогда
$$
y_t=\mu+u_t+z^{\gamma_n}_t {\xi}_t,\quad  u_t=\b u_{t-1} + {\e}_t, \quad  |\b|<1,\quad t \in \mathbb{Z,}
$$
наблюдаются $y_0,\ldots,y_n$. Тогда для $\hat\e_t$ из (1.6) имеем:
$$
\hat\e_t=\e_t-(\hat\b_n-\b) u_{t-1}-(1-\hat\b_n)(\hat \mu_n-\mu)+\alpha_t-(\hat\b_n-\b) z_{t-1}^{\gamma_n} \xi_{t-1},\quad t=1,\ldots,n,\eqno(2.1)
$$
где ради краткости мы положили $\alpha_t:=z^{\gamma_n}_t {\xi}_t-\b  z_{t-1}^{\gamma_n} \xi_{t-1}$.\\
В силу (2.1)
$$
\hat G_n(x)=\sum_{t-1}^n I(\e_t\leq x+(\hat\b_n-\b) u_{t-1}+(1-\hat\b_n)(\hat \mu_n-\mu)-\alpha_t+(\hat\b_n-\b) z_{t-1}^{\gamma_n} \xi_{t-1})\eqno(2.2)
$$
Свяжем с $\hat G_n(x)$ из (2.2)  еще одну о.э.ф.р.
$$
 G_n(x,\tau_1,\tau_2)=\sum_{t-1}^n I(\e_t\leq x+n^{-1/2}\tau_1 u_{t-1}+n^{-1/2}\tau_2-\alpha_t+n^{-1/2}\tau_1 z_{t-1}^{\gamma_n} \xi_{t-1}),
$$
и случайный процесс
$
u_n(x, \tau_1,\tau_2) :=$
$$
 n^{-1/2} \sum_{t=1}^n [I(\e_t \le x + n^{-1/2} \tau_1 u_{t-1}+n^{-1/2} \tau_2 - \alpha_t + n^{-1/2} \tau_1 z^{\gamma_n}_{t-1} \xi_{t-1}) -\\
 $$
 $$
  G( x + n^{-1/2} \tau_1 u_{t-1}+n^{-1/2} \tau_2 - \alpha_t + n^{-1/2} \tau_1 z^{\gamma_n}_{t-1} \xi_{t-1})].
$$
Пусть сигма-алгебра
$$\mathcal{F}_{t}: = \sigma\{\e_s,\, s \le t;\, \xi_i,\, z_i^{\gamma_n},\,i \le t+1\}.$$
Тогда слагаемы в $u_n(x, \tau_1,\tau_2)$ (при фиксированном $n$) образуют мартингал-разность относительно последовательности $\{\mathcal{F}_{t}\}$. \\

Далее $ 0 < \Theta < \infty$ фиксированно, $|\tau_i| \le \Theta$; $0\leq\Gamma<\infty,\,\Gamma$ фиксировано, $\gamma\leq\Gamma$. 
Утверждение теоремы 1.1 следует из следующих трех утверждений:
$$
\sup_{|\tau_i|\le\Theta}|u_n(x, \tau_1,\tau_2)-u_n(x, 0,0)|=o_P(1),\quad n\to \infty,\eqno (2.3)
$$
равномерно по $\gamma\le \Gamma$;
$$
\sup_{|\tau_i|\le\Theta}|n^{-1/2} \sum_{t=1}^n[G( x + n^{-1/2} \tau_1 u_{t-1}+n^{-1/2} \tau_2 - \alpha_t + n^{-1/2} \tau_1 z^{\gamma_n}_{t-1} \xi_{t-1})]-\eqno (2.4)
$$
$$
G( x  - \alpha_t )]-g(x)\tau_2|=o_P(1),\quad\to \infty,
$$
равномерно по $\gamma\le \Gamma$;

$$
n^{-1/2} \sum_{t=1}^n[I(\e_t\le x-\alpha_t)-I(\e_t\le x)]-\gamma \Delta(x,\Pi)|=o_P(1),\quad n\to \infty,\eqno (2.5)
$$
равномерно по $\gamma\le \Gamma$ c $\Delta(x,\Pi)=\sum_{j=0}^1[G(x+\b_j\xi_1)-G(x)] $ с $\b_0=-1,\b_1=\b$.\\
В самом деле, имеем тождество: $n^{1/2}[\hat G_n(x)-G_n(x)]=$
$$
u_n(x, n^{1/2}(\hat\b_n-\b),  n^{1/2}(1-\hat\b_n)(\hat \mu_n-\mu)-u_n(x,0, 0) +\eqno (2.6)
$$
$$
 n^{-1/2} \sum_{t=1}^n [  G( x + (\hat \b_n-\b) u_{t-1}+(1-\hat\b_n)(\hat \mu_n-\mu) - \alpha_t +(\hat \b_n-\b) z^{\gamma_n}_{t-1} \xi_{t-1})-G(x-\alpha_t)]+\eqno (2.7)
$$
$$
n^{-1/2} \sum_{t=1}^n[I(\e_t\le x-\alpha_t)-I(\e_t\le x)].\eqno (2.8)
$$
Пусть (2.3) верно. Тогда при любых положительных $\delta,\e$
$$
\sup_{\gamma\le \Gamma}\Prob(|u_n(x, n^{1/2}(\hat\b_n-\b),  n^{1/2}(1-\hat\b_n)(\hat \mu_n-\mu))-u_n(x,0, 0)|>\delta)\le 
$$
$$
\sup_{\gamma\le \Gamma}\Prob(|u_n(x, n^{1/2}(\hat\b_n-\b),  n^{1/2}(1-\hat\b_n)(\hat \mu_n-\mu))-u_n(x,0, 0)|>\delta,\, 
$$
$$|n^{1/2}(\hat \b_n-\b)|\le \Theta,\,|n^{1/2}(1-\hat\b_n)(\hat \mu_n-\mu)|\le \Theta)+
$$
$$
\sup_{\gamma\le \Gamma}\Prob((|n^{1/2}(\hat \b_n-\b)|> \Theta)+(|n^{1/2}(1-\hat\b_n)(\hat \mu_n-\mu)|> \Theta))\le
$$
$$
\sup_{\gamma\le \Gamma} \Prob(\sup_{|\tau_i|\le\Theta}|u_n(x, \tau_1,\tau_2)-u_n(x, 0,0)|>\delta)+
$$
$$
\sup_{\gamma\le \Gamma}\Prob((|n^{1/2}(\hat \b_n-\b)|> \Theta/2)+\sup_{\gamma\le \Gamma}\Prob(|n^{1/2}(1-\hat\b_n)(\hat \mu_n-\mu)|> \Theta/2)< \e
$$
при достаточно большом $\Theta>0$ для всех $n>n_0$.\\
Значит, (2.6) есть $o_P(1)$ равномерно по $\gamma\le \Gamma$.

Соверщенно аналогично показывается: соотношение (2.4) влечет, что (2.7) есть  $g(x)\tau_2+o_P(1)$  равномерно по $\gamma\le \Gamma$.
 
 Наконец, в \cite{Bold.Petr.} показано, что верно (2.5).\\
 Значит, (2.3)--(2.5) и (2.6)--(2.8) действительно влекут утверждение теоремы 1.1.
 
 Итак, будем доказывать (2.3)-(2.4). Начнем с (2.3) и дискретной аппроксимации $u_n(x, \tau_1,\tau_2)$. Напомним, $\Theta,\,\Gamma$ фиксированы, $|\tau_i|\le \Theta$, $\gamma\le \Gamma$.

Разобьём отрезок $[-\Theta n^{-1/2}, \Theta n^{-1/2}]$ на $3^{m_n}$ частей ($m_n$ - натуральные числа, $3^{m_n} \sim \ln n $ при $n \to \infty$) точками
$$
\eta_s =- \Theta n^{-1/2} + 2\Theta n^{-1/2} 3^{-m_n}s, s =0, 1,\dots,3^{m_n}.
$$
Пусть
\begin{align*}
\hat u_{ts}  = u_t[1 - 2 \Theta n^{-1/2} 3^{-m_n} \eta_s^{-1} I(u_t \le 0)], \\\
\tilde u_{ts}  = u_t[1 - 2 \Theta n^{-1/2} 3^{-m_n} \eta_s^{-1} I(u_t > 0)].
\end{align*}
Берем любые $\tau_1,\tau_2$ из отрезка $[-\Theta,\Theta]$. 
Выберем из точек $\{\eta_s\}$ точку $\eta_j$, ближайшую справа к $n^{-1/2} \tau_1$, тогда
$$0 \le \eta_j - n^{-1/2}\tau_1 \le 2 \Theta n^{-1/2} 3^{-m_n}.$$
Из введённых определений прямо следует:
\begin{gather}
\tag{2.9} 
\eta_j \tilde u_{t-1,j} \le  n^{-1/2} \tau_1 u_{t-1} \le \eta_j \hat u_{t-1,j}, \\
\tag{2.10} |\hat u_{ts}| \le 3|u_t|,\, |\tilde u_{ts}| \le 3|u_t|.
\end{gather}
Пусть $\eta_i$ будет ближайшая справа к $ n^{-1/2} \tau_2$ точка среди $\{\eta_s\}$. Тогда
$$
0 \le \eta_i-n^{-1/2} \tau_2\le 2\Theta n^{-1/2}3^{-m_n}.\eqno (2.11)
$$
Ради краткости положим еще
$$
\mu_t^{\pm} : = -\alpha_t \pm n^{-1/2}\Theta z^{\gamma_n}_{t-1} |\xi_{t-1}|.
$$

При $|\tau_i| \le \Theta$ в силу монотонности по $y$ $I(\e_t \le y)$ и $G(y)$ и силу (2.9) и (2.11) справедливы неравенства:
$$
u_n(x, \tau_1,\tau_2)-u_n(x, 0,0) \le 
$$
$$
 n^{-1/2} \sum_{t=1}^n [I(\e_t \le x + \eta_j \hat u_{t-1,j}+\eta_i +\mu_t^{+})-G( x + \eta_j \hat u_{t-1,j}+\eta_i +\mu_t^{+})  -
 $$
 $$
I(\e_t \le x -\alpha_t)+G( x  -\alpha_t)] + 
$$
$$
+n^{-1/2}\sum_{t=1}^n [G(x + \eta_j \hat u_{t-1,j}+\eta_i +\mu_t^{+})-G(x + \eta_j \tilde u_{t-1,j}+\eta_{i-1} +\mu_t^{-})],
$$
и, аналогично,
$$
u_n(x, \tau_1,\tau_2)-u_n(x, 0,0) \ge 
$$
$$
 n^{-1/2} \sum_{t=1}^n [I(\e_t \le x + \eta_j \tilde u_{t-1,j}+\eta_{i-1} +\mu_t^{-})-G(x + \eta_j \tilde u_{t-1,j}+\eta_{i-1} +\mu_t^{-})-
 $$
 $$
I(\e_t \le x -\alpha_t)+G( x  -\alpha_t)] -
$$
$$
n^{-1/2}\sum_{t=1}^n [G(x + \eta_j \hat u_{t-1,j}+\eta_i +\mu_t^{+})-G(x + \eta_j \tilde u_{t-1,j}+\eta_{i-1} +\mu_t^{-})].
$$
Из последних двух неравенств следует:
$$
\sup_{|\tau_i| \le \Theta}|u_n(x, \tau_1,\tau_2)-u_n(x, 0,0)| \le 
$$
$$
\max_{j,i}| n^{-1/2} \sum_{t=1}^n [I(\e_t \le x + \eta_j \hat u_{t-1,j}+\eta_i +\mu_t^{+})-G( x + \eta_j \hat u_{t-1,j}+\eta_i +\mu_t^{+})  -\eqno (2.12)
 $$
 $$
I(\e_t \le x -\alpha_t)+G( x  -\alpha_t)]| + 
$$
$$
\max_{j,i}|n^{-1/2} |\sum_{t=1}^n [I(\e_t \le x + \eta_j \tilde u_{t-1,j}+\eta_{i-1} +\mu_t^{-})-G(x + \eta_j \tilde u_{t-1,j}+\eta_{i-1} +\mu_t^{-})-\eqno(2.13)
 $$
 $$
I(\e_t \le x -\alpha_t)+G( x  -\alpha_t)]| +
$$
$$
\max_{j,i}n^{-1/2}\sum_{t=1}^n [G(x + \eta_j \hat u_{t-1,j}+\eta_i +\mu_t^{+})-G(x + \eta_j \tilde u_{t-1,j}+\eta_{i-1} +\mu_t^{-})].\eqno(2.14)
$$
Дискретная аппроксимаци завершена. 
Для доказательства соотношения (2.3) достаточно показать, что выражения в (2.12)-(2.14) есть $o_P(1),\,n\to\infty,$ равномерно по $\gamma\le \Gamma.$
\newtheorem{Lemma}{Лемма}[section]
\begin{Lemma}
{\sl
Выражения (2.12)--(2.13) есть $o_p(1)$ при  $n \to \infty$ равномерно по $\gamma \le \Gamma$.
}
\end{Lemma}
{\bf Доказательство.}
Покажем, что (2.12) есть $o_P(1)$ равномерно по $\gamma \le \Gamma$, для (2.13) рассуждения аналогичны. Положим
$$
\nu_t(j,i) := I(\e_t \le x + \eta_j \hat u_{t-1,j} +\eta_i+ \mu^{+}_t) - G(x + \eta_j \hat u_{t-1,j} +\eta_i+ \mu^{+}_t) - I(\e_t \le x -\alpha_t) + G(x -\alpha_t).
$$
Тогда (2.12) есть
$$
\max_{j,i}|\sum_{t=1}^n \nu_t(j,i)|.
$$
Пусть, как раньше, сигма-алгебра $\mathcal{F}_{t} = \sigma\{\e_s,\, s \le t;\, \xi_i,\,z_i^{\gamma_n},\, i \le t+1\}$.
Очевидно, последовтельность $\{\nu_t(j,i), \mathcal{F}_{t}\}, \,\,t=1,\ldots,n$, образует мартингал-разность. Значит, $\{\nu_t(j,i)\}$ последовательность центрированных и некоррелированных величин.

Далее воспользуемся известным неравенством: при $x_1, x_2 \in \mathbb{R}^1$
$$
 \E|I(\e_1 \le x_1) - G(x_1) - I(\e_1 \le x_2) + G(x_2)|^2 \le |G(x_1) - G(x_2)|.
$$
В силу этого неравенства, формулы Тейлора и   и (2.10)  
$$
\E\nu_t^2(j,i)= \E\E(\nu_t^2(j,i)/ \mathcal{F}_{t-1}) \le \E| G( x + \eta_j \hat u_{t-1,j} +\eta_i+ \mu^{+}_t) - G(x-\alpha_t)| \le
$$
$$
 \E| G( x +  \mu^{+}_t) - G(x-\alpha_t)| +
\sup_x g(x) \E|\eta_j \hat u_{t-1,j}+\eta_i|\le
$$
$$
\E| G( x +  \mu^{+}_t) - G(x-\alpha_t)|+c n^{-1/2} .
$$
Здесь через $c$ и далее через $c_1,c_2,\ldots$ обозначены константы, не зависящие от $n,j,i,t,\gamma$.\\
В силу формулы полной вероятности с гипотезами

$$
H_{00}=(z_0^{\gamma_n}=z_1^{\gamma_n}=0),\,H_{10}=(z_0^{\gamma_n}=1,\,z_1^{\gamma_n}=0),
$$
$$
H_{01}=(z_0^{\gamma_n}=0,\,z_1^{\gamma_n}=1),\,H_{11}=(z_0^{\gamma_n}=z_1^{\gamma_n}=1).
$$
получаем:
$$
 \E | G( x +  \mu^{+}_t)) - G(x-\alpha_t)|= \sum_{i,j} \E\{ [G(x + \mu_1^+) - G(x -\alpha_t)]/H_{ij}\}P(H_{ij})\le c_1 n^{-1/2}.
$$
Значит, $\E\nu_t^2(j,i)\le c_2 n^{-1/2}$, и 
$$
\E[n^{-1/2}\sum_{t=1}^n \nu_t(j,i)]^2=n^{-1}\sum_{t=1}^n \E\nu_t^2(j,i)\le c_2n^{-1/2}.
$$
 Отсюда, в силу нераенства Чебышева
$$
\Prob(\max_{j,i} |n^{-1/2}\sum_{t=1}^n \nu_t(j,i)| > \delta) \le \sum_{j,i} \Prob(|n^{-1/2}\sum_{t=1}^n \nu_t(j,i)| > \delta) \le
$$
$$
 \sum_{j,i} \delta^{-2}\E[n^{-1/2}\sum_{t=1}^n \nu_t(j,i)]^2  \le c_2 \delta^{-2}( 3^{m_n}+1) ^2n^{-1/2} = o(1)
$$
равномерно по $\gamma \le \Gamma$ в силу выбора последовательности $\{m_n\}$. Лемма 2.1 доказана.
\begin{Lemma}
{\sl
Выражение в (2.14) есть $o_p(1)$ при $n \to \infty$ равномерно по $\gamma \le \Gamma$.
}
\end{Lemma}
Доказательство леммы 2.2 лишь несущественными техническими деталями отличается от доказательства леммы 3.2 в \cite{Bold.Petr.}, 
а потому опущено.

Итак, соотношение (2.3) доказано.
\begin{Lemma}
{\sl
Утверждение (2.4) справедливо.
}
\end{Lemma}
{\bf Доказательство. }
При $|\tau_i|\le \Theta$ в силу монотонности $G(y)$ справедливы неравенства:
$$
n^{-1/2} \sum_{t=1}^n[G( x + n^{-1/2} \tau_1 u_{t-1}+n^{-1/2} \tau_2 - \alpha_t + n^{-1/2} \tau_1 z^{\gamma_n}_{t-1} \xi_{t-1})-
G( x  - \alpha_t )]\le
$$
$$
n^{-1/2} \sum_{t=1}^n[G( x + n^{-1/2} \tau_1 u_{t-1}+n^{-1/2} \tau_2  + \mu_t^+)-G( x  - \alpha_t )],\eqno(2.15)
$$
$$
n^{-1/2} \sum_{t=1}^n[G( x + n^{-1/2} \tau_1 u_{t-1}+n^{-1/2} \tau_2 - \alpha_t + n^{-1/2} \tau_1 z^{\gamma_n}_{t-1} \xi_{t-1})-
G( x  - \alpha_t )]\ge
$$
$$
n^{-1/2} \sum_{t=1}^n[G( x + n^{-1/2} \tau_1 u_{t-1}+n^{-1/2} \tau_2 + \mu_t^-)-G( x  - \alpha_t )].\eqno(2.16)
$$
Для доказательства леммы достаточно показать, что (2.15) и (2.16) равнмерно по  $\gamma \le \Gamma,\,|\tau_i| \le \Theta$ есть $g(x)\tau_2+o_P(1)$. 
Докажем это для (2.15), для (2.16) рассуждение аналогично. Итак, в силу формулы Тейлора (2.15) есть
$$
n^{-1/2} \sum_{t=1}^n[G( x +\mu_t^+)-G( x  - \alpha_t )]+\tau_1 n^{-1}\sum_{t=1}^ng( x +\mu_t^+)u_{t-1}+\tau_2 n^{-1}\sum_{t=1}^ng( x +\mu_t^+)
+o_P(1)\eqno(2.17)
$$
равнмерно по  $\gamma \le \Gamma,\,|\tau_i| \le \Theta$, т.к. вторые члены в разложении Тейлора есть, очевидно, $O(n^{-1/2})$ равномерно по $\gamma \le \Gamma,\,|\tau_i| \le \Theta$.

Первая сумма в (2.17) неотрицательна, ее среднее в силу формулы полной вероятности с гипотезами $H_{ij}$ издоказательства леммы 2.1 равно
$$
n^{1/2}\E [ G( x +  \mu^{+}_1)) - G(x-\alpha_1)]= \sum_{i,j=0}^1 n^{1/2} \E\{ [G(x + \mu_1^+) - G(x -\alpha_1)]/H_{ij}\}P(H_{ij})=
$$
$$
n^{1/2} \gamma_n (1-\gamma_n)\E [ G( x +  \b \xi_0+n^{-1/2} \Theta |\xi_0|)- G(x+\b \xi_0)] +\eqno(2.18)
$$
$$
n^{1/2} \gamma_n ^2\E [ G( x-\xi_1 +  \b \xi_0+n^{-1/2} \Theta |\xi_0|)- G(x-\xi_1+\b \xi_0)]=o(1),\quad \eqno(2.19)
$$
равномерно по $\gamma \le \Gamma$. Действительно, (2.18) есть $o(1)$ в силу теоремы Лебега, а (2.19) есть $O(n^{-1/2})$ равномерно по $\gamma \le \Gamma$.
 Значит, в силу неравенства Чебышева 1-ое слагаемое в (2.17) есть $o_P(1)$ равномерно по $\gamma \le \Gamma$ при 
 $n\to \infty$.

Вторая сумма в (2.17) при $|\tau_1| \le \Theta$ не больше 
$$
\Theta n^{-1}\sum_{t=1}^n |g(x+\mu_t^+)-g(x)||u_{t-1}|+\Theta g(x)|n^{-1}\sum_{t=1}^n u_{t-1}|=o_P(1) \eqno(2.20)
$$ 

равномерно по  $\gamma \le \Gamma,\,|\tau_1| \le \Theta$, т.к. среднее значение первого слагаемого в (2.20) есть $O(n^{-1/2})$ в силу формулы полной вероятности, а во втором слагаемом  $n^{-1}\sum_{t=1}^n u_{t-1}=o_P(1)$.

Наконец, третья сумма в (2.17) есть 
$$
\tau_2 n^{-1}\sum_{t=1}^n [g( x +\mu_t^+)-g(x)]+\tau_2 g( x )=\tau_2 g( x )+o_P(1),\quad n\to \infty,\eqno(2.21)
$$
равнoмерно по  $\gamma \le \Gamma,\,|\tau_2| \le \Theta$, т.к. среднее от супремума по $|\tau_2|\le \Theta$ модуля 1-го слагаемого в (2.21) есть $O(n^{-1/2})$ в силу формулы полной вероятности.

Лемма 2.3 доказана.

Итак, соотношения (2.3)--(2.5) доказаны, а потому теорема 1.1 полностью доказана.
\section{Доказательство теоремы 1.2.}
Доказательство теоремы 1.2 проводится по той же схеме, что и теоремы 1.1. Дальше $p=1$, $\{\e_t\}$ имеют ф.р. $A_n(x)$ из (1.7), $\hat G_n(x),\,G(x,\tau_1,\tau_2)$ топределяются так же, как в Разделе 2. 
Введм случайный процесс
$
u_n^A(x, \tau_1,\tau_2) :=$
$$
 n^{-1/2} \sum_{t=1}^n [I(\e_t \le x + n^{-1/2} \tau_1 u_{t-1}+n^{-1/2} \tau_2 - \alpha_t + n^{-1/2} \tau_1 z^{\gamma_n}_{t-1} \xi_{t-1}) -\\
 $$
 $$
 A_n( x + n^{-1/2} \tau_1 u_{t-1}+n^{-1/2} \tau_2 - \alpha_t + n^{-1/2} \tau_1 z^{\gamma_n}_{t-1} \xi_{t-1})].
$$
Аналогично доказательству теоремы 1.1 утверждение теоремы 1.2 следует из следующих трех утверждений:
$$
\sup_{|\tau_i|\le\Theta}|u_n^A(x, \tau_1,\tau_2)-u_n^A(x, 0,0)|=o_P(1),\quad n\to \infty,\eqno (3.1)
$$
равномерно по $\gamma\le \Gamma$;
$$
\sup_{|\tau_i|\le\Theta}|n^{-1/2} \sum_{t=1}^n[A_n( x + n^{-1/2} \tau_1 u_{t-1}+n^{-1/2} \tau_2 - \alpha_t + n^{-1/2} \tau_1 z^{\gamma_n}_{t-1} \xi_{t-1})-
$$
$$
A_n( x  - \alpha_t )]-g_0(x)\tau_2|=o_P(1),\quad n\to \infty,\eqno (3.2)
$$
равномерно по $\gamma\le \Gamma$;

$$
n^{-1/2} \sum_{t=1}^n[I(\e_t\le x-\alpha_t)-I(\e_t\le x)]-\gamma \Delta_0(x,\Pi)|=o_P(1),\quad n\to \infty,\eqno (3.3)
$$
равномерно по $\gamma\le \Gamma$ c $\Delta_0(x,\Pi)=\sum_{j=0}^1[G_0(x+\b_j\xi_1)-G_0(x)] $ и $\b_0=-1,\b_1=\b$.\\
Соотношение (3.1) доказывается так же, как соотношение (2.3), надо лишь в доказательстве (2.3) везде заменить ф.р. $G(x)$ на ф.р. $A_n(x)$. Разумеется, 
используется тот факт, что при Условии (А) ф.р. $A_n(x)$ дважды дифференцируема и обе производные ограничены по $n,x$. Кроме того, $\sup_n \E\e_1^2<\infty$, откуда 
следует для строго стационарног решения уравнения авторегрессии: $\sup_n \E u_1^2<\infty$.\\
Справедливость соотношения (3.3) установлена в \cite{Bold.2019}. Докажем (3.2).

\begin{Lemma}
{\sl
Утверждение (3.2) справедливо.
}
\end{Lemma}
{\bf Доказательство. }
При $|\tau_i|\le \Theta$ в силу моотонности $A_n(y)$ справедливы неравенства:
$$
n^{-1/2} \sum_{t=1}^n[A_n( x + n^{-1/2} \tau_1 u_{t-1}+n^{-1/2} \tau_2 - \alpha_t + n^{-1/2} \tau_1 z^{\gamma_n}_{t-1} \xi_{t-1})-
A_n( x  - \alpha_t )]\le
$$
$$
n^{-1/2} \sum_{t=1}^n[A_n( x + n^{-1/2} \tau_1 u_{t-1}+n^{-1/2} \tau_2  + \mu_t^+)-A_n( x  - \alpha_t )],\eqno(3.4)
$$
$$
n^{-1/2} \sum_{t=1}^n[A_n( x + n^{-1/2} \tau_1 u_{t-1}+n^{-1/2} \tau_2 - \alpha_t + n^{-1/2} \tau_1 z^{\gamma_n}_{t-1} \xi_{t-1})-
A_n( x  - \alpha_t )]\ge
$$
$$
n^{-1/2} \sum_{t=1}^n[A_n( x + n^{-1/2} \tau_1 u_{t-1}+n^{-1/2} \tau_2  + \mu_t^-)-A_n( x  - \alpha_t )].\eqno(3.5)
$$
Для доказательства леммы достаточно показать, что выражения в (3.4) и (3.5) равномерно по  $\gamma \le \Gamma,\,|\tau_i| \le \Theta$ есть $g_0(x)\tau_2+o_P(1)$. 
Докажем это для (3.4), для (3.5) рассуждение аналогично. Итак, в силу формулы Тейлора (3.4) есть
$$
n^{-1/2} \sum_{t=1}^n[A_n( x +\mu_t^+)-A_n( x  - \alpha_t )]+\tau_1 n^{-1}\sum_{t=1}^n a_n( x +\mu_t^+)u_{t-1}+\tau_2 n^{-1}\sum_{t=1}^n a_n( x +\mu_t^+)
+o_P(1)\eqno(3.6)
$$
равнмерно по  $\gamma \le \Gamma,\,|\tau_i| \le \Theta$. Здесь $a_n(x)=A_n'(x)$. 

Первая сумма в (3.6) есть $o_P(1)$ равномерно по  $\gamma \le \Gamma$, это показывается совершенно аналогично доказательству соотношения 
$$
n^{-1/2} \sum_{t=1}^n[G( x +\mu_t^+)-G( x  - \alpha_t )]=o_P(1),\quad n\to \infty,
$$ 
установленного в доказательстве леммы 2.3.
 
Вторая сумма в (3.6) не больше
$$
\Theta n^{-1}\sum_{t=1}^n |a_n(x+\mu_t^+)-a_n(x)||u_{t-1}|+\Theta a_n(x)|n^{-1}\sum_{t=1}^n u_{t-1}|=o_P(1) \eqno(3.7)
$$ 
равномерно по  $\gamma \le \Gamma,\,|\tau_1| \le \Theta$, т.к. среднее значение первого слагаемого в (3.7) есть $O(n^{-1/2})$, а во втором слагаемом  $n^{-1}\sum_{t=1}^n u_{t-1}=o_P(1),\,$.

Наконец, третья сумма в (3.6) есть 
$$
\tau_2 n^{-1}\sum_{t=1}^n [a_n( x +\mu_t^+)-a_n(x)]+\tau_2 a_n( x )=\tau_2 g_0( x )+o_P(1),\quad n\to \infty,\eqno(3.8)
$$
равнoмерно по  $\gamma \le \Gamma,\,|\tau_2| \le \Theta$, т.к. среднее от супремума по $|\tau_2|\le \Theta$ модуля первого слагаемого в (3.8) есть $O(n^{-1/2})$ в силу формулы полной вероятности, а $a_n(x)=g_0(x)+o(1),\,n\to\infty$.

Лемма 3.1 доказана.

Соотношения (3.1)--(3.3) проверены. Теорема 1.2 доказана.


\end{document}